\documentstyle[twoside,psfig]{article}

\oddsidemargin=\evensidemargin 
\catcode`\@=11 \long\def\@makefntext#1{ \protect\noindent \hbox to
3.2pt {\hskip-.9pt
$^{{\eightrm\@thefnmark}}$\hfil}#1\hfill}       %CAN BE USED

\def\ps@myheadings{\let\@mkboth\@gobbletwo      %SIZE OF R/H NOS.
\def\@oddhead{\hbox{}
\rightmark\hfil\eightrm\thepage}
\def\@oddfoot{}\def\@evenhead{\eightrm\thepage\hfil
\leftmark\hbox{}}\def\@evenfoot{}
\def\sectionmark##1{}\def\subsectionmark##1{}}

\catcode`@=11                        %SIZE OF OPENING PAGE NO.
\def\ps@plain{\let\@mkboth\@gobbletwo
     \def\@oddhead{}\def\@oddfoot{\eightrm\hfil\thepage
     \hfil}\def\@evenhead{}\let\@evenfoot\@oddfoot}

\newcounter{sectionc}\newcounter{subsectionc}\newcounter{subsubsectionc}
\renewcommand{\section}[1] {\vspace{12pt}\addtocounter{sectionc}{1}
\setcounter{subsectionc}{0}\setcounter{subsubsectionc}{0}\noindent
    {\tenbf\thesectionc. #1}\par\vspace{5pt}}
\renewcommand{\subsection}[1] {\vspace{12pt}\addtocounter{subsectionc}{1}
    \setcounter{subsubsectionc}{0}\noindent
    {\bf\thesectionc.\thesubsectionc.
    {\kern1pt \bfit #1}}\par\vspace{5pt}}
\renewcommand{\subsubsection}[1] {\vspace{12pt}
    \addtocounter{subsubsectionc}{1}
    \noindent
    {\tenrm\thesectionc.\thesubsectionc.\thesubsubsectionc. {\kern1pt
    \it #1}}\par\vspace{5pt}}
\newcommand{\nonumsection}[1] {\vspace{12pt}\noindent{\tenbf #1}
    \par\vspace{5pt}}

\newcounter{appendixc}
\newcounter{subappendixc}[appendixc]
\newcounter{subsubappendixc}[subappendixc]

\renewcommand{\appendix}[1] {\vspace{12pt}  %Appendix A. Heading +
    \refstepcounter{appendixc}      %Appendix A (11/5/94)
    \setcounter{figure}{0}
    \setcounter{table}{0}
    \setcounter{lemma}{0}
    \setcounter{theorem}{0}
    \setcounter{corollary}{0}
    \setcounter{definition}{0}
    \setcounter{equation}{0}
    \renewcommand{\thefigure}{\Alph{appendixc}.\arabic{figure}}
    \renewcommand{\thetable}{\Alph{appendixc}.\arabic{table}}
    \renewcommand{\theappendixc}{\Alph{appendixc}}
    \renewcommand{\thelemma}{\Alph{appendixc}.\arabic{lemma}}
    \renewcommand{\thetheorem}{\Alph{appendixc}.\arabic{theorem}}
    \renewcommand{\thedefinition}{\Alph{appendixc}.\arabic{definition}}
    \renewcommand{\thecorollary}{\Alph{appendixc}.\arabic{corollary}}
    \renewcommand{\theequation}{\Alph{appendixc}.\arabic{equation}}
    \noindent{\tenbf Appendix \theappendixc #1}\par\vspace{5pt}}

\newcommand{\textlineskip}{\baselineskip=13pt}
\newcommand{\smalllineskip}{\baselineskip=10pt}

\def\abstracts#1#2#3#4{{
    \centering{\begin{minipage}{4.5in}\footnotesize\baselineskip=10pt
    \centerline{ABSTRACT}
    \parindent=15pt #1\par
    \parindent=15pt #2\par
    \parindent=15pt #3\par
    \parindent=15pt #4\par
    \end{minipage}}\par}}

\renewenvironment{thebibliography}[1]
    {\frenchspacing
     \ninerm\baselineskip=11pt
     \begin{list}{[\arabic{enumi}]}
    {\usecounter{enumi}\setlength{\parsep}{0pt}
     \setlength{\leftmargin 13.7pt}{\rightmargin 0pt} %[1--9] ITEMS
     \setlength{\itemsep}{0pt} \settowidth
    {\labelwidth}{[#1]}\sloppy}}{\end{list}}

\newcounter{itemlistc}
\newcounter{romanlistc}
\newcounter{alphlistc}
\newcounter{arabiclistc}

\newcommand{\fcaption}[1]{
        \refstepcounter{figure}
        \setbox\@tempboxa = \hbox{\footnotesize Fig.~\thefigure. #1}
        \ifdim \wd\@tempboxa > 5in
           {\begin{center}
        \parbox{5in}{\footnotesize\smalllineskip Fig.~\thefigure. #1}
            \end{center}}
        \else
             {\begin{center}
             {\footnotesize Fig.~\thefigure. #1}
              \end{center}}
        \fi}

\newcommand{\tcaption}[1]{
        \refstepcounter{table}
        \setbox\@tempboxa = \hbox{\footnotesize Table~\thetable. #1}
        \ifdim \wd\@tempboxa > 5in
           {\begin{center}
        \parbox{5in}{\footnotesize\smalllineskip Table~\thetable. #1}
            \end{center}}
        \else
             {\begin{center}
             {\footnotesize Table~\thetable. #1}
              \end{center}}
        \fi}

\def\pmb#1{\setbox0=\hbox{#1}
    \kern-.025em\copy0\kern-\wd0
    \kern.05em\copy0\kern-\wd0
    \kern-.025em\raise.0433em\box0}

\def\fnt#1#2{\footnotetext{\kern-.3em
    {$^{\mbox{\scriptsize #1}}$}{#2}}}

\def\fpage#1{\begingroup
\voffset=.3in
\thispagestyle{empty}\begin{table}[b]\centerline{\footnotesize #1}
    \end{table}\endgroup}

\def\runninghead#1#2{\pagestyle{myheadings}
\markboth{{\protect\footnotesize\it{\quad #1}}\hfill}
{\hfill{\protect\footnotesize\it{#2\quad}}}}

\font\tenrm=cmr10  \font\tenbf=cmbx10
\font\bfit=cmbxti10 at 10pt \font\ninerm=cmr9 \font\nineit=cmti9
 \font\eightrm=cmr8

\newtheorem{theorem}{Theorem}   %P.S. -- CONNECTED TO PROP, REMARK ETC...
\newtheorem{definition}{Definition}

\def\@begintheorem#1#2{\trivlist    %6/9/94
    \item[\hskip\labelsep{\bf #1\ #2.}]}
\def\@opargbegintheorem#1#2#3{\trivlist
    \item[\hskip\labelsep{\bf #1\ #2\ (#3).}]}

        {\setcounter{itemlistc}{0}      %\NOINDENT\BULLET
     \begin{list}{$\bullet$}        %AUTOMATICALLY
    {\usecounter{itemlistc}         %UNLESS OTHERWISE SPECIFIED
     \leftmargin10pt           %IE.  0PT = MOVE OUTSIDE LEFT MARGIN
     \setlength{\parsep}{0pt}
     \setlength{\itemsep}{0pt}     %ADD EXTRA VSPACE BETWEEN ITEMS
    }}{\end{list}}
    {\setcounter{romanlistc}{0}     %\NOINDENT (i)
     \begin{list}{$($\roman{romanlistc}$)$} %AUTOMATICALLY
    {\usecounter{romanlistc}        %UNLESS OTHERWISE SPECIFIED
     \leftmargin18pt %(ii)      %IE. 21PT = MOVE MORE TO THE RIGHT-->
     \setlength{\parsep}{0pt}
     \setlength{\itemsep}{0pt}  %ADD EXTRA VSPACE BETWEEN ITEMS
     \settowidth{\labelwidth}{#1}
    }}{\end{list}}
    {\setcounter{enumii}{0}         %\NOINDENT (a)
     \begin{list}{$($\alph{enumii}$)$}  %AUTOMATICALLY
    {\usecounter{enumii}            %UNLESS OTHERWISE SPECIFIED
     \leftmargin18pt        %IE.  0PT = MOVE OUTSIDE LEFT MARGIN
     \setlength{\parsep}{0pt}
     \setlength{\itemsep}{0pt}  %ADD EXTRA VSPACE BETWEEN ITEMS
     \settowidth{\labelwidth}{#1}
    }}{\end{list}}
\def\qed{\hbox{${\vcenter{\vbox{            %HOLLOW SQUARE
   \hrule height 0.4pt\hbox{\vrule width 0.4pt height 6pt
   \kern5pt\vrule width 0.4pt}\hrule height 0.4pt}}}$}}

\def\theequation{\thesectionc.\arabic{equation}}  %FOR SETTING EQ.~(1.1)
\begin{document}
\setlength{\textheight}{7.7truein}  %for 2nd page onwards

\runninghead{Braid Family Representatives} {Braid Family
Representatives}

\normalsize\textlineskip \thispagestyle{empty}
\setcounter{page}{1}

%\copyrightheading{}         %{Vol.~0, No.~0 (1999) 00--00}

%\vspace*{0.88truein}

\fpage{1} \centerline{\bf BRAID FAMILY REPRESENTATIVES}
\vspace*{0.37truein}
\centerline{\footnotesize SLAVIK JABLAN and RADMILA SAZDANOVI\'
C$^\dag $} \baselineskip=12pt \centerline{\footnotesize\it The
Mathematical Institute, Knez Mihailova
35,}\centerline{\footnotesize\it P.O.Box 367, 11001
Belgrade,}\centerline{\footnotesize\it Serbia $\&$ Montenegro}
\centerline{\footnotesize\it jablans@mi.sanu.ac.yu}
\centerline{\footnotesize\it seasmile@galeb.etf.bg.ac.yu$^\dag $}

%\vspace*{0.225truein} \publisher{}

\vspace*{0.21truein} \abstracts{After defining reduced minimum
braid word and criteria for a braid family representative,
different braid family representatives are derived, and a
correspondence between them and families of knots and links given
in Conway notation is established.}{}{}{}

%\vspace*{10pt}
%\keywords{The contents of the keywords}

%\textlineskip          %) USE THIS MEASUREMENT WHEN THERE IS
%\vspace*{12pt}         %) NO SECTION HEADING

\vspace*{1pt}\textlineskip  %) USE THIS MEASUREMENT WHEN THERE IS
\section{Introduction}    %) A SECTION HEADING
\vspace*{-0.5pt}

In the present article  Conway notation [1,2,3,4] will be used
without any additional explanation. A {\it braid-modified Conway
notation} is introduced in Section 1,  for a better understanding
of the correspondence between braid family representatives
($BFR$s) and families of knots and links ($KL$s) given in Conway
notation.

Minimum braids are defined, described, generated and presented in
tables for knots up to ten crossings and oriented links up to nine
crossings by T.~Gittings [5]. T.~Gittings used them for studying
graph trees, amphicheirality, unknotting numbers and periodic
tables of $KL$s.

Since knots are 1-component links, the term $KL$ will be used for
both knots and links.

In Section 2 we define a {\it reduced braid word}, describe
general form for all reduced braid words with $s=2$ strands,
generate all braid family representatives of two-strand braids,
and establish a correspondence between them and families of $KL$s
given in Conway notation. In Section 3 we consider the same
problem for $s\ge3$. In Section 4 some applications of minimum
braids [5] and braid family representatives are discussed. All
computations are made using the knot-theory program {\it LinKnot}
written by the authors [6], the extension of the program {\it
Knot2000} by M.~Ochiai and N.~Imafuji [7].

\vspace*{1pt}\textlineskip  %) USE THIS MEASUREMENT WHEN THERE IS
\section{Reduced Braid Words and Minimum Families of Braids with $s=2$}    %) A SECTION HEADING
\vspace*{-0.5pt}

We use the standard definition of a braid and description of
minimum braids given by T.~Gittings [5]. Instead of $a\ldots a$,
where a capital or lower case letter $a$ appears $p$ times, we
write $a^p$; $p$ is the degree of $a$ ($p\in N$). It is also
possible to work with negative powers, satisfying the
relationships: $A^{-p}=a^p$, $a^{-p}=A^p$. A number of strands is
denoted by $s$, and a length of a braid word by $l$.

The operation $a^2\rightarrow a$ applied on any capital or lower
case letter $a$ is called {\it idempotency}. To every braid word
we can apply the operation of idempotency until a reduced braid
word is obtained.

\medskip

\begin{definition}
A {\it reduced braid word} is a braid word with degree of every
capital or lower case letter equal to 1.
\end{definition}

\medskip

By an opposite procedure, {\it braid word extension}, from every
reduced braid word we obtain all braid words that can be derived
from it by assigning a degree (that can be greater then 1) to
every letter. In this case, a reduced braid word plays a role of a
{\it generating braid word}.

A braid word with one or more parameters denoting degrees greater
then one represents a {\it family of braid words}. If values of
all parameters are equal 2, it will be called a {\it source
braid}.

For the {\it minimality of reduced braids} we are using the
following criteria:

\begin{enumerate}
\item minimum number of braid crossings;
\item minimum number of braid strands;
\item minimum binary code for {\it alternating} braid crossings.
\end{enumerate}

According to the first and second criterion minimal reduced braids
are the shortest reduced braids with a smallest as possible number
of different letters among all equivalent reduced braids
representing certain $KL$. A binary code for any braid crossing
can be generated by assigning a zero for an alternating, and a one
for a non-alternating crossing. Hence, a priority will be given to
alternating braids, and then to braids that differ from them as
low as possible. Analogous minimality criteria can be applied to
source braids.

\medskip

\begin{definition}
Among the set of all braid families representing the same $KL$
family, the {\it braid family representative} ($MFB$) is the one
that has the following properties:

\begin{enumerate}
\item minimum number of braid crossings;
\item minimum {\it reduced braid};
\item minimum {\it source braid}.
\end{enumerate}
\end{definition}

These criteria are listed in descending order of importance for
determining $BFR$s.

Our definition of $BFR$s results in some fundamental differences
with regard to minimum braids, defined by T.~Gittings [5]. Some
members of $BFR$s will be minimum braids, but not necessarily.

For example, the minimum braid of the link $.2\,1:2$ ($9_{11}^3$
in Rolfsen [4]) is 9:03-05a $AAbACbACb$ [5,Table 2]. According to
the second $BFR$ criterion it will be derived from the generating
minimum braid $AbAbACbC$ corresponding to the link $.2\,1$
($8_{13}^2$), and not from the non-minimum generating braid
$AbACbACb$ corresponding to the same link. Hence, to the
three-component link $.2\,1:2$ ($9_{11}^3$) obtained as the first
member of $BFR$ $AbA^pbACbC$ for $p=2$ will correspond the braid
$AbAAbACbC$, that is not a minimum braid according to the minimum
braid criteria [5].

The third criterion: {\it minimum source braid} enables us to
obtain $KL$s of a certain family from a single $BFR$, and not from
several different $BFR$s. For example, applying this criterion,
$KL$s $.3.2.2\, 0$, $.2.3.2\,0$ and $.2.2.3\,0$ belonging to the
same $KL$ family $.r.p.q\,0$ will be obtained from the single
$BFR$ $A^pbA^qbAb^r$. Otherwise, using the minimum braid criteria
[5], the knot $.3.2.2\,0$ will be obtained from the family
$A^pbAb^qAb^r$, three-component link $.2.3.2\,0$ will be obtained
from $A^pbA^qbAb^r$, and the knot $.2.2.3\,0$ will be obtained
from $A^pbA^qb^rAb$ for $p=3$, $q=2$, $r=2$. Source braids
corresponding to the families $A^pbAb^qAb^r$, $A^pbA^qbAb^r$ and
$A^pbA^qb^rAb$ are $A^2bAb^2Ab^2$, $A^2bA^2bAb^2$ and
$A^2bA^2b^2Ab$, respectively, and the second source braid is
minimal. Hence, the representative of the $KL$ family $.r.p.q\,0$
is $BFR$ $A^pbA^qbAb^r$.

According to this, to every $BFR$ can be associated a single
corresponding family of $KL$s given in Conway notation and {\it
vice versa}.

An overlapping of $KL$ families obtained from $BFR$s can occur
only at their beginnings. For example, distinct $BFR$s
$AbA^pbACbC$ and $A^pbCbAbCb$, giving $KL$ families $.2\,1:p$ and
$.p\,1:2$, respectively, for $p=2$ will have as a joint member
aforementioned three-component link $.2\,1:2$ ($9_{11}^3$).
According to the second $BFR$ criterion, it will be derived from
the minimum generating braid $AbAbACbC$, and not from $AbACbACb$.
Hence, $BFR$ $AbA^pbACbC$ giving $KL$s of the form $.2\,1:p$
begins for $p=2$, and $A^pbCbAbCb$ giving $KL$s of the form for
$.p\,1:2$ begins for $p=3$. In this way, all ambiguous cases can
be solved.

Every $KL$ is algebraic (if its basic polyhedron is $1^*$) or
polyhedral, so according to this criterion, all $KL$s are divided
into two main categories: algebraic and polyhedral. Since to every
member of a $BFR$ corresponds a single $KL$, we can introduce the
following definition:

\medskip

\begin{definition}
An alternating $BFR$ is {\it polyhedral} {\it iff} its
corresponding $KL$s are polyhedral. Otherwise, it is {\it
algebraic}. A non-alternating $BFR$ is {\it polyhedral} {\it iff}
its corresponding alternating $BFR$ is polyhedral. Otherwise, it
will be called {\it algebraic}.
\end{definition}

\medskip

The division of non-alternating $BFR$s into algebraic and
polyhedral does not coincide with the division of the
corresponding $KL$s [1,2,3], because minimum number of braid
crossings is used as the first criterion for the $BFR$s. Accepting
{\it minimum reduced braid universe} [5] as the first criterion,
all $KL$s derived from the basic polyhedron $.1$ will be
algebraic, because they can be represented by non-alternating
minimal (but not minimum [5]) algebraic braids. E.g., the
alternating knot $.2.2\,0$ ($8_{16}$) with the polyhedral braid
$A^2bA^2bAb$ can be represented as the algebraic knot
$(-3,2)\,(3,-2)$ with the corresponding algebraic braid
$A^3b^2a^2B^3$. In this case, to the knot $8_{16}$ corresponds
algebraic braid $A^3b^2a^2B^3$ that reduces to $AbaB$, and not
$A^2bA^2bAb$ that reduces to $AbAbAb$.

Another solution of this discrepancy is changing the definition of
an algebraic $KL$ into the following:

\begin{definition}
$KL$ is algebraic if it has an algebraic {\it minimum crossing
number} representation.
\end{definition}

\noindent In this case, all $KL$s derived from the basic
polyhedron $.1$ (with Conway symbols beginning with a dot) will be
polyhedral $KL$s, because their minimum crossing number
representations are polyhedral.

We will consider only $BFR$s corresponding to prime $KL$s.

It is easy to conclude that every 1-strand $BFR$ is of the form
$A^p$, with the corresponding $KL$ family $p$ in Conway notation.

\bigskip

\begin{theorem}
Every reduced $BFR$ with $s=2$ is of the form $(Ab)^n$, $n\ge 2$.
\end{theorem}

\bigskip

This $BFR$ corresponds to the knot $2\,2$ and to the family of
basic polyhedra $.1=6^*$, $8^*$, $10^*$, $12^*$ (or 12A according
to A.~Caudron [3]), {\it etc}. For $n\ge 3$ all of them are
$n$-antiprisms. Let us notice that the first member of this
family, the knot $2\,2$, is not an exception: it is an antiprism
with two digonal bases.

\bigskip

\begin{theorem}
All algebraic alternating $KL$s with $s=2$ are the members of the
following families:

\medskip

$p\,1\,2 $ with the $BFR$ $A^pbAb$ ($p\ge 1$);

$p\,1\,1\,q$ with the $BFR$ $A^pbAb^q$ ($p\ge q\ge 2$);

$p,q,2$ with the $BFR$ $A^pbA^qb$ ($p\ge q\ge 2$);

$p,q,r\,1$ with the $BFR$ $A^pbA^qb^r$ ($r\ge 2, p\ge q \ge 2$);

$(p,r)\,(q,s)$ with the $BFR$ $A^pb^qA^rb^s$

($p,q,r,s\ge 2, p\ge r, p\ge s, s\ge q$ and if $p=s$, then $r\ge
q$).
\end{theorem}

\bigskip

Minimum braids include one additional braid ($A^pb^qAb^r$) in the
case of algebraic alternating $KL$s with $s=2$.

Alternating polyhedral $KL$s with $s=2$ are given in the following
table, each with its $BFR$. $KL$s in this table are given in
"standard" Conway notation (that is "standardized" for knots with
$n\le 10$ and links with $n\le 9$ crossings according to Rolfsen's
book [4]). This table can be extended to an infinite list of
antiprismatic basic polyhedra $(2n)^*$ described by the $BFR$s
$(Ab)^n$, $n\ge 3$ and $BFR$s with $s=2$ obtained as their
extensions.

\bigskip

{\bf Table 1}

\medskip

{\bf Basic polyhedron $.1=6^*$}

\medskip

\halign{ # \hfill &  #  &  #  & # &  #  & $\quad \quad \quad \quad
$ #  & # \hfill & # & # &  #  & #  \cr
 $A^pbAbAb$   & & $.p$ & & (1) & & $A^pbAbA^qb^r$& &$r:p\,0:q\,0$ & (7) & \cr
 $A^pbAbAb^q$ &  &$.p.q$ &  & (2) & &$A^pbAb^qA^rb^s$ & & $p.s.r.q$& (8) & \cr
 $A^pbA^qbAb$ &  &$.p.q\,0$ & & (3) & & $A^pbA^qbA^rb^s$& & $q\,0.p.r\,0.s\,0$ & (9) &\cr
 $A^pbAb^qAb$ &  & $.p:q\,0$ & & (4) & & $A^pbA^qb^rAb^s$ & & $.p.s.r\,0.q\,0$ & (10) & \cr
 $A^pbA^qbAb^r$  & &$.r.p.q\,0$ & & (5) & &   $A^pbA^qb^rA^sb^t$ & & $p.t.s.r.q$  & (11) &\cr
 $A^pbA^qbA^rb$ & & $p:q:r$ & & (6) & &   $A^pb^qA^rb^sA^tb^u$ & &$p.q.r.s.t.u$ & (12) & \cr}

\bigskip

If we apply minimum braid criteria [5], we need to add ten braids
for the basic polyhedron $.1=6^*$: (1') $A^pbAb^qAb^r$, (2')
$A^pbA^qb^rAb$, (3') $A^pb^qAbAb^r$, (4') $A^pbA^qb^rA^sb$, (5')
$A^pb^qAbA^rb^s$, (6') $A^pb^qAb^rAb^s$, (7') $A^pb^qA^rbAb^s$,
(8') $A^pb^qAb^rA^sb^t$, (9') $A^pb^qA^rbA^sb^t$, (10')
$A^pb^qA^rb^sAb^t$. Applying $BFR$ criteria, according to the
minimum source braid criterion all $KL$s obtained from the braids
(1') and (2') will be obtained from $BFR$ (5), $KL$s obtained from
(3') will be obtained from (7), $KL$s obtained from (4') and (6')
will be obtained from (9), $KL$s obtained from (5') and (7') will
be obtained from (8), and $KL$s obtained from (8'), (9') and (10')
will be obtained from (11). Using minimum braid criteria [5], we
need to make analogous additions to all classes of $BFR$s
considered in this paper.

\bigskip

For the basic polyhedron $8^*$ we have:

\bigskip

{\bf Basic polyhedron $8^*$}

\medskip

\halign{ # \hfill &  #  &  #  & $\quad \quad $ #  & # \hfill & # &
# \cr
 $A^pbAbAbAb$ & & $8^*p$ & & $A^pbA^qbAb^rAb^s$  & &$8^*p:q:.r:s$  \cr
 $A^pbAbAbAb^q$ & & $8^*p.q$ & &$A^pbAb^qA^rbAb^s$& &$8^*p.s:.r.q$  \cr
 $A^pbA^qbAbAb$ & & $8^*p:q$ & & $A^pbA^qbA^rbA^sb$& &$8^*p:s:r:q$  \cr
 $A^pbAbAb^qAb$ & & $8^*p:.q$ & & $A^pbAbA^qb^rA^sb^t$& & $8^*p.t.s.r.q$ \cr
 $A^pbAbA^qbAb$ & & $8^*p::q$ & & $A^pbA^qbAb^rA^sb^t$& &$8^*p.t.s.r:.q$   \cr
 $A^pbA^qbAbAb^r$ & & $8^*p.r::.q$ & & $A^pbA^qb^rA^sbAb^t$& &$8^*p:q.r.s:.t$ \cr
 $A^pbAbA^qbAb^r$ & & $8^*p.r:.q$ & & $A^pbA^qbA^rbA^sb^t$& & $8^*p.t.s:r:q$ \cr
 $A^pbA^qbA^rbAb$ & & $8^*p:q:r$ & & $A^pbA^qbA^rb^sAb^t$ & &$8^*p.t:s.r:q$   \cr
 $A^pbA^qbAb^rAb$ & & $8^*p:.r:.q$ & &$A^pbAb^qA^rb^sA^tb^u$& &$8^*p.u.t.s.r.q$ \cr
 $A^pbAbAbA^qb^r$ & & $8^*p.r.q$ & & $A^pbA^qbA^rb^sA^tb^u$& & $8^*p.u.t.s.r:q$ \cr
 $A^pbAbAb^qA^rb^s$ & & $8^*p.s.r.q$ & & $A^pbA^qb^rA^sbA^tb^u$& & $8^*p.u.t:s.r.q$ \cr
 $A^pbA^qbAbA^rb^s$ & & $8^*p.s.r::q$ & &$A^pbA^qb^rA^sb^tAb^u$ & & $8^*p:q.r.s.t:u$ \cr
 $A^pbAb^qA^rb^sAb$ & & $8^*p:.s.r.q$ & & $A^pbA^qb^rA^sb^tA^ub^v$ & & $8^*p.v.u.t.s.r.q$ \cr
 $A^pbA^qb^rAbAb^s$ & & $8^*p.s::r.q$ & & $A^pb^qA^rb^sA^tb^uA^vb^w$ & & $8^*p.q.r.s.t.u.v.w$ \cr
 $A^pbA^qbA^rbAb^s$ & & $8^*p.s:.r:q$  & &  & &  \cr}

 \bigskip

Trying to better understand the correspondence between $BFR$s and
Conway symbols of $KL$s, we can introduce {\it modified Conway
notation}. Most of $KL$s can be given in Conway notation by
several different symbols (and this is the main disadvantage of
Conway notation). In a similar way as with the classical notation,
where every $KL$ is given by its place in knot tables, we need to
use some "standard" code, according to the notation introduced in
the original Conway's paper [2] and in the papers or books
following it [1,3,4]. For example, the same polyhedral knot $.p$
can be given by $..p$, $:p$, $:.p$, $\ldots $, or even as $6^*p$,
$6^*.p$, $6^*:.p$, $\ldots $

Working with $BFR$s we introduce a {\it braid-modified Conway
notation} that will be more suitable for denoting $KL$s obtained
from $BFR$s. We are trying to have a same degree $p$ at the first
position of a braid, and as the first element of Conway symbol
corresponding to it. Whenever possible, the order of degrees will
be preserved in the corresponding Conway symbol. By using this
notation, we can recognize a very simple pattern for $BFR$s
derived from the generating minimum braids of the form $(Ab)^n$:
by denoting in a Conway symbol corresponding to a given braid
every sequence of single letters of a length $k$ by $k+1$ dots, we
obtain the Conway symbol of a given braid. In order to recognize
this pattern for $KL$s derived from basic polyhedra, first we need
to use only one basic polyhedron $6^*$ with $n=6$ crossings, and
not two of them ($.1$ and $6^*$). In this case, the Table 1 will
look as follows:

\bigskip

{\bf Basic polyhedron $6^*$}

\medskip

\halign{ # \hfill &  #  &  #  & $\quad \quad \quad \quad $ #  & #
\hfill & # & # \cr
 $A^pbAbAb$   & & $6^*p$ & & $A^pbAbA^qb^r$& &$6^*p::q.r$ \cr
 $A^pbAbAb^q$ &  &$6^*p::.q$ &  &$A^pbAb^qA^rb^s$ & & $6^*p:.q.r.s$\cr
 $A^pbA^qbAb$ &  &$6^*p:q$ & & $A^pbA^qbA^rb^s$& & $6^*p:q:r.s$\cr
 $A^pbAb^qAb$ &  & $6^*p:.q$ & & $A^pbA^qb^rAb^s$ & & $6^*p:q.r:s$ \cr
 $A^pbA^qbAb^r$  & &$6^*p:q:.r$ & &   $A^pbA^qb^rA^sb^t$ & & $6^*p:q.r.s.t$  \cr
 $A^pbA^qbA^rb$ & & $6^*p:q:r$ & &   $A^pb^qA^rb^sA^tb^u$ & &$6^*p.q.r.s.t.u$ \cr}

\bigskip

\noindent and for the basic polyhedron $8^*$ we have:

\bigskip

{\bf Basic polyhedron $8^*$}

\medskip

\halign{ # \hfill &  #  &  #  & $\quad \quad $ #  & # \hfill & # &
# \cr
 $A^pbAbAbAb$ & & $8^*p$ & & $A^pbA^qbAb^rAb^s$  & &$8^*p:q:.r:s$  \cr
 $A^pbAbAbAb^q$ & & $8^*p:::.q$ & &$A^pbAb^qA^rbAb^s$& &$8^*p:.q.r:.s$  \cr
 $A^pbA^qbAbAb$ & & $8^*p:q$ & & $A^pbA^qbA^rbA^sb$& &$8^*p:q:r:s$  \cr
 $A^pbAbAb^qAb$ & & $8^*p::.q$ & & $A^pbAbA^qb^rA^sb^t$& & $8^*p::q.r.s.t$ \cr
 $A^pbAbA^qbAb$ & & $8^*p::q$ & & $A^pbA^qbAb^rA^sb^t$& &$8^*p:q:.r.s.t$   \cr
 $A^pbA^qbAbAb^r$ & & $8^*p:q::.r$ & & $A^pbA^qb^rA^sbAb^t$& &$8^*p:q.r.s:.t$ \cr
 $A^pbAbA^qbAb^r$ & & $8^*p::q:.r$ & & $A^pbA^qbA^rbA^sb^t$& & $8^*p:q:r:s.t$ \cr
 $A^pbA^qbA^rbAb$ & & $8^*p:q:r$ & & $A^pbA^qbA^rb^sAb^t$ & &$8^*p:q:r.s:t$   \cr
 $A^pbA^qbAb^rAb$ & & $8^*p:q:.r$ & &$A^pbAb^qA^rb^sA^tb^u$& &$8^*p:.q.r.s.t.u$ \cr
 $A^pbAbAbA^qb^r$ & & $8^*p:::q.r$ & & $A^pbA^qbA^rb^sA^tb^u$& & $8^*p:q:r.s.t.u$ \cr
 $A^pbAbAb^qA^rb^s$ & & $8^*p::.q.r.s$ & & $A^pbA^qb^rA^sbA^tb^u$& & $8^*p:q.r.s:t.u$ \cr
 $A^pbA^qbAbA^rb^s$ & & $8^*p:q::r.s$ & &$A^pbA^qb^rA^sb^tAb^u$ & & $8^*p:q.r.s.t:u$ \cr
 $A^pbAb^qA^rb^sAb$ & & $8^*p:.q.r.s$ & & $A^pbA^qb^rA^sb^tA^ub^v$ & & $8^*p:q.r.s.t.u.v$ \cr
 $A^pbA^qb^rAbAb^s$ & & $8^*p:q.r::s$ & & $A^pb^qA^rb^sA^tb^uA^vb^w$ & & $8^*p.q.r.s.t.u.v.w$ \cr
 $A^pbA^qbA^rbAb^s$ & & $8^*p:q:r:.s$  & &  & &  \cr}

\medskip

Unfortunately, it is not possible to express every family of $KL$s
in the braid-modified Conway notation. Another problem is that it
strongly differs from the standard Conway notation. Therefore, the
braid-modified Conway notation is used only when after some slight
modification standard Conway symbols remained completely
understandable to a reader familiar with them.

In the same way, it is possible to continue with the derivation of
$BFR$s from basic polyhedra with a higher number of crossings.

\medskip

Hence, we conclude that:

\medskip

{\bf Corollary} All alternating $KL$s with $s=2$ are described by
Theorem 2 and by an infinite extension of Table 1.

\medskip

From alternating $BFR$s we obtain non-alternating $BFR$s by
crossing changes. This way, from $BFR$s derived from the
generating minimum braid $(Ab)^2$ we obtain the following families
of non-alternating $BFR$s and corresponding new $KL$ families:

\medskip

\halign{ # \hfill &  #  &  #  & $\quad \quad \quad \quad $ #  & #
\hfill & # & # \cr
 $A^pBaB$ & & $(p-1)\,3$ & & $A^pba^qb^r$ & & $p,(q-1)\,1,-(r+1)$  \cr
  $A^pba^qb$ & & $p,(q-1)\,1,-2$ & &$A^pBA^qB^r$ &  & $p,q,-r\,1$ \cr
  $A^pBA^qB$& & $p,q,-2$  &  &$A^pB^qa^rB^s$ & & $(-p,r)\,(q,s)$ \cr
  $A^pBaB^q$ & & $(p-1)\,2\,q$ & &$A^pB^qA^rB^s$  & &$(p,r)-(q,s)$ \cr
 $A^pB^qaB^r$ & & $p-1,q,r+$ & &  & & \cr}

 \medskip

In the same way, we can derive non-alternating $BFR$s with $s=2$
from the generating $BFR$ $(Ab)^n$, $n\ge 3$.

\vspace*{1pt}\textlineskip  %) USE THIS MEASUREMENT WHEN THERE IS
\section{Braid Family Representatives with $s\ge 3$}    %) A SECTION HEADING
\vspace*{-0.5pt}

In order to continue derivation of $BFR$s and corresponding $KL$s
for $s\ge 3$ first we derive  all different reduced minimum braid
words. It is possible to establish general construction rules for
generating minimum braid words.

\medskip

\begin{definition}
For a given generating minimum braid word $W=wL$ that ends with a
capital or lower case letter $L$, a replacement of $L$ by a word
$w_1$ in $W$ will be called {\it extending by replacement}. An
addition of the word $w_1$ to $W$ is {\it extending by addition}.
The both operations are {\it extending operations}.
\end{definition}

\medskip

\begin{definition}
Let $W=wL_s$ and $w_1=L_{s+1}L_sL_{s+1}$ be generating minimum
braids with $s$ and $s+1$ strings, where $L_s$ denotes $s$th
letter and $L_{s+1}$ denotes $(s+1)$th letter. The word extending
operations obtained this way will be called, respectively, {\it
$(s+1)$-extending by replacement}, and {\it $(s+1)$-extending by
addition}. The both operations are {\it $(s+1)$-extending
operations}.
\end{definition}

\medskip

For example, the first operation applied on $AbAb$ gives $AbACbC$,
and the other $AbAbCbC$.

The $(s+1)$-extending by replacement is sufficient for
construction of generating minimum braids for a given $s$, with
$l=2s$, corresponding to $KL$s of the form $2\,\ldots\,2=2^s$,
where 2 occurs $s$ times. For $2\le s\le 6$ as the result we
obtain: $AbAb$, $AbACbC$, $AbACbdCd$, $AbACbdCEdE$, $AbACbdCEdfEf$
$\ldots$

The generating minimum braids for given $s$, with $l=3s-2$,
corresponding to $KL$s of the form
$2\,1\,\ldots\,1\,2=2\,1^{3s-6}\,2$, where 1 occurs $3s-6$ times,
can be obtained using only $(s+1)$-extension by addition. For
$3\le s\le 6$ we obtain: $AbAbCbC$, $AbAbCbCdCd$, $AbAbCbCdCdEdE$,
$AbAbCbCdCdEdEfEf$ $\ldots $

Applying the same procedure, from $A^3$ we obtain the series
$A^3BaB$, $A^3BaBCbC$, $A^3BaBCbCDcD$, $A^3BaBCbCDcDEdE$ $\ldots
$, corresponding to the knots $3\,2$, $5\,2$, $7\,2$, $9\,2$
$\ldots$

Analogously, starting with $w_1=AbAbCbdCd$ and using the
$(s+1)$-extension by replacement, the generating minimum braids
with $l=2s+1$, corresponding to $KL$s of the form $2\,2\,1\,\ldots
1\,2=2^2\,1^{2s-5}\,2$ are obtained for given $s$.

However, in order to exhaust all possibilities, all combinations
of $(s+1)$-extending operations are used for derivation of reduced
minimum braids.

\medskip

\begin{theorem}
Every generating algebraic minimum braid can be derived from
$AbAb$ by a recursive application of $(s+1)$-extending operations.
\end{theorem}

\medskip

The minimal generating braid words for $s\le 5$ with their
corresponding $KL$s are given in the following table:

\medskip

\halign{ # \hfill &  #  &  # \hfill  & # & # \hfill & # & # \hfill
\cr
 $s=1$   &  & $l=1$ &  & $A$ &  & 1 \cr
    &  &  &  &  &  &  \cr
 $s=2$   &  & $l=4$ &  & $AbAb$ &  & $2\,2$ \cr
    &  &  &  &  &  &  \cr
 $s=3$   &  & $l=6$ &  & $AbACbC$ &  & $2\,2\,2$ \cr
 $s=3$   &  & $l=7$ &  & $AbAbCbC$ &  & $2\,1\,1\,1\,2$ \cr
    &  &  &  &  &  &  \cr
 $s=4$   &  & $l=8$ &  & $AbACbdCd$ &  & $2\,2\,2\,2$  \cr
 $s=4$   &  & $l=9$ &  & $AbAbCbdCd$ &  & $2\,2\,1\,1\,1\,2$ \cr
 $s=4$   &  & $l=10$ &  & $AbAbCbCdCd$ &  & $2\,1\,1\,1\,1\,1\,1\,2$ \cr
    &  &  &  &  &  &  \cr
 $s=5$   &  & $l=10$ &  & $AbACbdCEdE$ &  & $2\,2\,2\,2\,2$ \cr
 $s=5$   &  & $l=11$ &  & $AbAbCbdCdEdE$ &  & $2\,2\,2\,1\,1\,1\,2$ \cr
 $s=5$   &  & $l=11$ &  & $AbACbCdCEdE$ &  & $2\,2\,1\,1\,1\,2\,2$ \cr
 $s=5$   &  & $l=12$ &  & $AbAbCbCdCEdE$ &  & $2\,2\,1\,1\,1\,1\,1\,1\,2$ \cr
 $s=5$   &  & $l=12$ &  & $AbAbCbdCdEdE$ &  & $2\,1\,1\,1\,2\,1\,1\,1\,2$ \cr
 $s=5$   &  & $l=13$ &  & $AbAbCbCdCdEcE$ &  & $2\,1\,1\,1\,1\,1\,1\,1\,1\,1\,2$ \cr}

\medskip

In the case of polyhedral generating minimum braid words it is
also possible to make generalizations. We have already considered
the infinite class of generating polyhedral minimum braid words
$(Ab)^n$ with $s=2$. The first infinite class with $s=3$ will be
$(Ab)^{n-1}ACbC$, with the corresponding $KL$s of the form
$(2n)^*2\,1\,0$.

Every $BFR$ can be derived from a generating minimum braid by
assigning a degree (that can be greater then 1) to every letter.

For $s=3$ there are two generating alternating algebraic minimum
braid words:

$AbACbC$, $l=6$, with the corresponding link $2\,2\,2$;

$AbAbCbC$, $l=7$, with the corresponding knot $2\,1\,1\,1\,2$,

\noindent that generate prime $KL$s.

From $AbACbC$ we derived 17 alternating $BFR$s and their
corresponding families of $KL$s, given in the following table:

\medskip

\halign{ # \hfill &  #  &  #  & $ \quad \quad \quad $ #  & #
\hfill & # & # \cr
 $A^pbACbC$ & & $p\,1\,2\,2$& & $A^pbA^qCb^rC$& & $(p,q)\,(r,2+)$  \cr
 $AbACb^pC$ & & $p,2,2+$    & & $A^pb^qA^rCb^sC$& &$(p,r)\,(q,2,s)$   \cr
 $A^pbACb^qC$ & &$p\,1,q,2+$& & $A^pb^qACb^rC^s$& &$p\,1,q,s\,1,r$   \cr
 $A^pbA^qCbC$ & &$p,q,2\,2$ & & $A^pbA^qCb^rC^s$& & $(p,q)\,(r,s\,1+)$  \cr
 $A^pbACbC^q$ & &$p\,1\,2\,1\,q$ & & $A^pbA^qC^rbC^s$& & $(p,q)\,2\,(r,s)$  \cr
 $Ab^pACb^qC$ & &$p,2,q,2$& & $A^pb^qA^rCb^sC^t$& &$(p,r)\,(q,t\,1,s)$   \cr
  $A^pb^qACb^rC$ & &$p\,1,q,r,2$& & $A^pbA^qC^rb^sC^t$& &  $(p,q),s,(t,r)+$ \cr
 $A^pbACb^qC^r$ & &$p\,1,q,r\,1+$& & $A^pb^qA^rC^sb^tC^u$& & $(p,r),q,(u,s),t$  \cr
  $A^pbA^qCbC^r$ & &$p,q,r\,1\,2$& &   \cr}

 \medskip

The next generating alternating algebraic minimum braid $AbAbCbC$
of the length 7, with $s=3$, gives the following results:
\medskip

\halign{ # \hfill &  #  &  #  & $\quad \quad $ #  & # \hfill & # &
# \cr
 $A^pbAbCbC$ & & $p\,1\,1\,1\,1\,2$& &$A^pb^qAbCb^rC^s$ & &
 $(p\,1,q)\,1\,(s\,1,r)$\cr
 $AbAbCb^pC$ & & $p,2\,1\,1,2$& &$A^pbA^qb^rCb^sC$ & &
 $(p,q)\,1\,r\,(2,s)$\cr
 $AbAb^pCbC$ & & $2\,1\,p\,1\,2$& & $A^pbA^qb^rCbC^s$& & $p,q,s\,1\,1\,r\,1$  \cr
 $A^pb^qAbCbC$ & & $p\,1,q,2\,1\,1$& & $A^pbA^qbC^rbC^s$& & $(p,q)\,1\,1\,1\,(r,s)$\cr
 $A^pbA^qbCbC$ & & $p,q,2\,1\,1\,1$& &  $A^pbA^qbCb^rC^s$& & $(p,q)\,1\,1\,(r,s\,1)$\cr
 $A^pbAb^qCbC$ & & $p\,1\,1\,q\,1\,2$& & $A^pbAb^qCb^rC^s$& & $p\,1\,1\,q,r,s\,1$  \cr
 $A^pbAbCbC^q$ & & $p\,1\,1\,1\,1\,1\,q$& & $A^pbAbC^qb^rC^s$& & $(p\,1\,1\,1,r)\,(q,s)$\cr
  $AbAb^pCb^qC$ & & $2\,1\,p,q,2$& & $A^pb^qA^rbCb^sC$& &$(p,r),q,(2,s)\,1$ \cr
$Ab^pAbCb^qC$ & & $(p,2)\,1\,(q,2)$& & $A^pb^qAb^rCb^sC$&
&$(p\,1,q)\,r\,(2,s)$ \cr
 $A^pb^qA^rbCbC$ & & $(p,r)\,(q,2\,1\,1)$& &$A^pb^qA^rbCb^sC^t$& & $(p,r),q,(t\,1,s)\,1$  \cr
 $A^pb^qAb^rCbC$ & & $p\,1,q,2\,1\,r$& & $A^pb^qAb^rCb^sC^t$& & $(p\,1,q)\,r\,(t\,1,s)$  \cr
 $A^pb^qAbCb^rC$ & & $(p\,1,q)\,1\,(r,2)$& & $A^pbA^qb^rC^sbC^t$& & $(p,q)\,1\,r\,1\,(s,t)$  \cr
 $A^pbA^qbCb^rC$ & & $(p,q)\,1\,1\,(2,r)$& & $A^pbA^qb^rCb^sC^t$& & $(p,q)\,1\,r\,(s,t\,1)$  \cr
 $A^pbA^qb^rCbC$ & & $p,q,2\,1\,r\,1$& & $A^pbA^qbC^rb^sC^t$& & $(p,q)\,1\,1,(t,r),s$  \cr
 $A^pbA^qbCbC^r$ & & $p,q,r\,1\,1\,1\,1$& & $A^pbAb^qC^rb^sC^t$& & $(p\,1\,1\,q,s)\,(r,t)$  \cr
 $A^pbAb^qCb^rC$ & & $p\,1\,1\,q,r,2$& & $A^pb^qA^rb^sCb^tC$& & $(p,r),q,(t,2)\,s$  \cr
 $A^pbAb^qCbC^r$ & & $p\,1\,1\,q\,1\,1\,r$& & $A^pb^qA^rb^sCb^tC^u$& & $((p,r),q)\,s\,(u\,1,t)$  \cr
 $A^pbAbCb^qC^r$ & & $p\,1\,1\,1,q,r\,1$& & $A^pb^qA^rbC^sb^tC^u$& & $((p,r),q)\,1\,((u,s),t)$  \cr
 $Ab^pAb^qCb^rC$ & & $(p,2)\,q\,(r,2)$& & $A^pbA^qb^rC^sb^tC^u$& & $(p,q)\,1\,r\,((u,s),t)$  \cr
 $A^pb^qA^rb^sCbC$ & &
 $(p,r)\,(q,2\,1\,s)$& & $A^pb^qA^rb^sC^tb^uC^v$& & $((p,r),q)\,s\,((u,t),v)$  \cr}

 \medskip

Except $AbACbC$ and $AbAbCbC$, all generating minimum braids with
$s=3$ are polyhedral.

For $s=3$ and $l\le 12$, the polyhedral generating braids and
their corresponding $KL$s are given in the following table, with
the notation for basic polyhedra with 12 crossings according to
A.~Caudron [3]:

\footnotesize

\medskip

\halign{ # \hfill &  #  &  # \hfill  & # & # \hfill & # & $\quad
$# & # \hfill &  #  &  # \hfill  & # & # \hfill & # \cr
 $l=8$  &  & $AbAbACbC$ &  & $.2\,1$  &  & $l=11$ &  & $AbAbAbAbCbC$ &  & $8^*2\,1\,1$ \cr
 $l=8$  &  & $AbCbAbCb$ &  & $.2:2$ &  &  $l=11$ &  & $AbAbAbCbCbC$ &  & $11^{***}$ \cr
  &  &  &  &  &  & $l=11$ &  & $AbAbACbAbCb$ &  & $10^{**}.2\,0$ \cr
 $l=9$ &  & $AbAbCbAbC$ &  & $8^*2\,0$ &  & $l=11$ &  & $AbAbACbACbC$ &  & $11^{**}$ \cr
 $l=9$ &  & $AbAbAbCbC$ &  & $.2\,1\,1$ &  & $l=11$ & & $AbAbCbACbCb$& &$11^*$ \cr
 $l=9$ &  & $AbACbACbC$ &  & $9^*$ &  &  &  & &  & \cr
  &  &  &  &  &  &  $l=12$ &  & $AbAbAbAbACbC$ &  &  $10^*2\,1\,0$ \cr
 $l=10$ &  & $AbAbAbCbCb$ &  & $.2\,1\,2$  &  & $l=12$ &  & $AbAbAbACbCbC$ &  & $12$I \cr
 $l=10$ &  & $AbAbAbACbC$ &  & $8^*2\,1\,0$  &  & $l=12$ &  & $AbAbAbCbACbC$ &  & $12$F \cr
 $l=10$ &  & $AbAbACbAbC$ &  & $9^*.2$  &  & $l=12$ &  & $AbAbACbAbCbC$ &  & $12$H  \cr
 $l=10$ &  & $AbAbCbAbCb$ &  & $9^*2$  &  &  $l=12$ &  & $AbAbCbAbACbC$   &  &  $12$G\cr
 $l=10$ &  & $AbAbACbCbC$ &  & $10^{***}$ &  &  $l=12$  &  & $AbAbCbAbCbCb$&   &  $12$D\cr
 $l=10$  &  & $AbAbCbACbC$ &  & $10^{**}$   &  &$l=12$ &  & $AbCbAbCbAbCb$  &  & $12$C  \cr}

\medskip

\normalsize

From them, $BFR$s without duplications are derived. E.g., for
$l=8$, the generating minimum braid $.2\,1$ gives 70 $BFR$s, and
$.2:2$ gives 19 $BFR$s. Overlapping of those families can occur
only if all parameters are equal 2, i.e., for source braids and
source $KL$s corresponding to them. According to the minimality
criteria, all those source braids will belong to the first $BFR$.
The generating minimum braid $.2\,1$ gives the following $BFR$s:

\footnotesize

\medskip

\halign{ # \hfill &  #  &  # \hfill  & # $\quad $  & # \hfill & #
& # \hfill \cr
 $A^pbAbACbC$ & & $.2\,1.p\,0$ & & $AbA^pbAC^qbC^r$ & & $.p:(q,r) 1$\cr
 $AbA^pbACbC$ & & $.2\,1:p$ & & $AbA^pb^qACb^rC$ & & $.p.q.(2,r)$\cr
 $AbAbACb^pC$ & & $.(p,2)$ & & $AbA^pb^qACbC^r$ & & $.r\,1\,\,1.q.p$\cr
 $AbAbACbC^p$ & & $.p\,1\,1$ & & $AbAb^pACb^qC^r$ & & $.p.(r\,1,q)$\cr
 $AbAb^pACbC$ & & $.2\,1.p$ & & $AbAb^pAC^qbC^r$ & & $.p.(r,q)\,1$\cr
 $A^pbA^qbACbC$ & & $.2\,1.p\,0.q$ & & $Ab^pA^qb^rACbC$ & & $2\,1.p.r\,0.q\,0$\cr
 $A^pbAbA^qCbC$ & & $2\,1\,0:p\,0:q\,0$& &  $Ab^pAb^qACb^rC$& & $p:q:(2,r)\,0$ \cr
 $A^pbAbACb^qC$ & & $.(2,q).p\,0$ & & $Ab^pAb^qACbC^r$ & & $q:p:r\,1\,1\,0$\cr
 $A^pbAbACbC^q$ & & $.q\,1\,1.p\,0$ & & $A^pbA^qbA^rCb^sC$ & &$(2,s).p\,0.r.q\,0$ \cr
 $A^pbAb^qACbC$ & & $.q.2\,1.p\,0$& & $A^pbA^qbA^rCbC^s$ & & $s\,1\,1\,0:r\,0.q.p\,0$ \cr
 $A^pb^qAbACbC$ & & $2\,1:p:q\,0$  & & $A^pbA^qbACb^rC^s$ & &$.(s\,1,r).p\,0.q$ \cr
 $AbA^pbACb^qC$ & & $.(2,q):p$& & $A^pbA^qbAC^rbC^s$ & & $(s,r)\,1:p.q\,0$\cr
 $AbA^pbACbC^q$ & & $.q\,1\,1:p$ & & $A^pbA^qb^rA^sCbC$ & & $q.r.s.2\,1\,0.p$\cr
 $AbA^pb^qACbC$ & & $.2\,1.q.p$  & & $A^pbA^qb^rACb^sC$ & &  $.(s,2).r.q.p\,0$\cr
 $AbAbACb^pC^q$ & & $.(q\,1,p)$ & & $A^pbA^qb^rACbC^s$ & & $.s\,1\,1.r.q.p\,0$\cr
 $AbAbAC^pbC^q$ & & $.(q,p)\,1$ & & $A^pbAbA^qCb^rC^s$ & & $(r,s\,1)\,0:q\,0:p\,0$\cr
 $AbAb^pACb^qC$ & & $.(q,2).p$  & & $A^pbAbA^qC^rbC^s$ & & $(s,r)\,1\,0:q\,0:p\,0$\cr
 $AbAb^pACbC^q$ & & $.q\,1\,1.p$ & & $A^pbAb^qA^rCb^sC$ & & $q.(2,s).r\,0.p$\cr
 $Ab^pAb^qACbC$ & & $p:q:2\,1\,0$ & & $A^pbAb^qA^rCbC^s$ & & $q.s\,1\,1.r\,0.p$\cr
 $A^pbA^qbA^rCbC$ & & $.2\,1.p\,0.q:r$ & & $A^pbAb^qACb^rC^s$ & & $.q.(s\,1,r).p\,0$\cr
 $A^pbA^qbACb^rC$ & & $.q.p\,0.(r,2)$ & & $A^pbAb^qAC^rbC^s$ & & $.q.(s,r)\,1.p\,0$\cr
 $A^pbA^qbACbC^r$ & & $.r\,1\,1.p\,0.q$ & & $A^pb^qA^rbACb^sC$ & & $r.q.p.(2,s)\,0$\cr
 $A^pbA^qb^rACbC$ & & $.2\,1.r.q.p\,0$ & & $A^pb^qA^rbACbC^s$ & & $r.q.p.s\,1\,1\,0$\cr
 $A^pbAbA^qCb^rC$ & & $p\,0:q\,0:(r,2)\,0$ & & $A^pb^qA^rb^sACbC$ & & $s.r.q.p.2 1 0$\cr
 $A^pbAbA^qCbC^r$ & & $r\,1\,1\,0:q\,0:p\,0$ & & $A^pb^qAbACb^rC^s$ & & $(s\,1,r):p:q\,0$\cr
 $A^pbAbACb^qC^r$ & & $.p.(r\,1,q)\,0$ & & $A^pb^qAbAC^rbC^s$ & & $(s,r)\,1:p:q\,0$\cr
 $A^pbAbAC^qbC^r$ & & $.p.(r,q)\,1\,0$ & & $A^pb^qAb^rA^sCbC$ & & $s.r.q.p.2\,1\,0$\cr
 $A^pbAb^qA^rCbC$ & & $q.2\,1.r\,0.p$ & & $A^pb^qAb^rACb^sC$ & & $p\,0.(s,2).q.r\,0$\cr
 $A^pbAb^qACb^rC$ & & $.q.(2,r).p\,0$ & & $A^pb^qAb^rACbC^s$ & & $p\,0.s\,1\,1.q.r\,0$\cr
 $A^pbAb^qACbC^r$ & & $.q.r\,1\,1.p\,0$ & & $AbA^pb^qACb^rC^s$ & & $.(r,s\,1).q.p$\cr
 $A^pb^qA^rbACbC$ & & $r.q.p.2\,1\,0$ & & $AbA^pb^qAC^rbC^s$ & & $.(s,r)\,1.q.p$\cr
 $A^pb^qAbACb^rC$ & & $p:(r,2):q\,0$ & & $Ab^pA^qb^rACb^sC$ & & $(s,2).p.r\,0.q\,0$\cr
 $A^pb^qAbACbC^r$ & & $r\,1\,1:p:q\,0$ & & $Ab^pA^qb^rACbC^s$ & & $s\,1\,1.r.p\,0.q\,0$\cr
 $A^pb^qAb^rACbC$ & & $p\,0.2\,1.q.r 0$ & & $Ab^pAb^qACb^rC^s$ & & $q:p:(s\,1,r)\,0$\cr
 $AbA^pbACb^qC^r$ & & $.p:(r\,1,q)$ & & $Ab^pAb^qAC^rbC^s$ & & $p:q:(r,s)\,1\,0$\cr}

\medskip

\normalsize

The generating minimum braid $.2:2$ gives the following $BFR$s:

\footnotesize

\medskip

\halign{ # \hfill &  #  &  # \hfill  & # $\quad $  & # \hfill & #
& # \hfill \cr
   $A^pbCbAbCb$ & & $.p\,1:2$ & & $A^pb^qCbAbC^rb$ & & $.r\,1.q\,0.p\,1$ \cr
   $Ab^pCbAbCb$ & & $.2.p\,0.2$ & & $A^pb^qCb^rAbCb$ & & $2.q\,0.r.p\,1\,0$ \cr
   $A^pb^qCbAbCb$ & & $.p\,1.q\,0.2$ & & $A^pbCbA^qbC^rb$ & & $.r\,1:(p,q)$ \cr
   $A^pbCbA^qbCb$ & & $.(p,q):2$ & & $Ab^pCbAb^qCb^r$ & & $2.p\,0.r.2\,0.q$ \cr
   $A^pbCbAbC^qb$ & & $.p\,1:q\,1$ & & $A^pb^qCbA^rb^sCb$ & & $.(p,r).s\,0.2.q\,0$ \cr
   $Ab^pCbAb^qCb$ & & $.p.2\,0.q.2\,0$ & & $A^pb^qCbA^rbC^sb$ & & $.s\,1.q\,0.(p,r)$ \cr
   $Ab^pCbAbCb^q$ & & $2.p\,0.q.2\,0$ & & $A^pb^qCbA^rbCb^s$ & & $(r,p).q\,0.s.2\,0$ \cr
   $A^pb^qCbA^rbCb$ & & $.2.q\,0.(p,r)$ & & $A^pb^qCbAb^rC^sb$ & & $.p\,1.r\,0.s\,1.q\,0$ \cr
   $A^pb^qCbAb^rCb$ & & $.p\,1.r\,0.2.q\,0$ & & $A^pb^qCbAbC^rb^s$ & & $p\,1.s\,0.q.r\,1\,0$ \cr
   $A^pb^qCbAbCb^r$ & & $p\,1.q\,0.r.2\,0$ & & $$ & & $$ \cr}

\medskip

\normalsize

In the same way, from all generating minimum braid words with
$s=3$, it is possible to derive alternating and non-alternating
$BFR$s and their corresponding families of $KL$s.

As the example, the following table contains non-alternating
$BFR$s with at most two parameters, derived from the minimum
reduced braid $AbACbC$:

\medskip

\halign{ # \hfill &  #  &  # \hfill  & # $\quad \quad \quad $  & #
\hfill & # & # \hfill \cr
  $A^pBacBc$ & & $(p-1)\,3\,2$ & & $A^pba^qCbC$ & & $p,2\,2,-q$ \cr
  $A^pBacBc^q$ & & $(p-1)\,3\,1\,q$ & & $A^pBA^qcBc$ & & $p,2\,1\,1,-(q-1)\,1$ \cr
  $A^pbACB^qC$ & & $p\,1,(q-1)\,1,2$ & & $Ab^pACB^qC$ & & $p,2,2,-q$ \cr
  $A^pBacB^qc$ & & $p-1,q,2++$ & & $Ab^pAcB^qc$ & & $p,2,-q,-2$ \cr
  $A^pbACB^qC$ & & $p\,1,(q-1)\,1,2$ & & $Ab^pAcb^qc$ & & $p,q,2,-2$ \cr
  $A^pbACbc^q$ & & $p\,1\,3\,(q-1)$ & & $AB^pACB^qC$ & & $p,q,-2,-2$ \cr
  $A^pBacBC^q$ & & $(p-1)\,4\,(q-1)$ & & $AB^pACb^qC$ & & $p,2,q,-2$ \cr}

\medskip

From the generating minimum braid word $W=(Ab)^n$ ($n\ge 2$), that
defines the family of basic polyhedra $(2n)^*$, by word extension
$w_1=CbACbC$ we obtain the second family of basic polyhedra $9^*$
($AbACbACbC$), $10^{**}$ ($AbAbCbACbC$), $11^{**}$
($AbAbACbACbC$), $12$F ($AbAbAbCbACbC$), {\it etc.}

The third family of basic polyhedra $10^{***}$ ($AbAbACbCbC$),
$11^{***}$ ($AbAbAbCbCbC$), $12$I ($AbAbAbACbCbC$), {\it etc.}, is
derived from $W=(Ab)^n$ ($n\ge 3$) for $w_1=CbCbC$.

In the same way, for $W=(Ab)^n$ ($n\ge 1$), $w_1=CbAbCbAbCb$ the
family of basic polyhedra beginning with $12$C ($AbCbAbCbAbCb$) is
obtained;

for $W=(Ab)^n$ ($n\ge 2$), $w_1=CbAbCbCb$ the family of basic
polyhedra beginning with $12$D ($AbAbCbAbCbCb$) is obtained;

for $W=(Ab)^n$ ($n\ge 2$), $w_1=CbAbACbC$  the family of basic
polyhedra beginning with $12$G ($AbAbCbAbACbC$) is obtained;

for $W=(Ab)^n$ ($n\ge 2$), $w_1=CbAbCbC$  the family of basic
polyhedra beginning with $12$H ($AbAbACbAbCbC$) is obtained, {\it
etc}.

\noindent Among them it is possible to distinguish subfamilies
obtained using extensions by replacing or by adding.

\medskip

\begin{theorem}
For $s=4$ generating algebraic minimum braids are:

 $AbACbdCd$, $l=8$, with the corresponding link $2\,2\,2\,2$,

 $AbAbCbdCd$, $l=9$, with the corresponding link
 $2\,2\,1\,1\,1\,2$,

 $AbAbCbCdCd$, $l=10$, with the corresponding knot
 $2\,1\,1\,1\,1\,1 \,1\,2$.

 \noindent All other generating minimum braid words with
 $s=4$ are polyhedral.
\end{theorem}

\medskip

For $s=4$ and $l\le 12$, the polyhedral generating braids and
their corresponding $KL$s are given in the following table, with
the notation for basic polyhedra with 12 crossings according to
A.~Caudron [3]:

\medskip

\footnotesize

\halign{ # \hfill &  #  &  # \hfill  & # & # \hfill & # & $\quad
 $# & # \hfill &  #  &  # \hfill  & # & # \hfill & #
\cr
  $l=10$  &  & $AbAbACbdCd$ &  & $.2\,2\,1$  &  & $l=12$ &  & $AbAbACbdCdCd$& & $12$J \cr
  $l=10$  &  & $AbACbCbdCd$ &  & $.2\,1.2\,1$  &  &   $l=12$ &  & $AbAbACdCbCdC$  & & $11^{***}:.2\,0$  \cr
  $l=10$  &  & $AbACbdCbdC$ &  & $.2\,1:2\,1\,0$  &  &  $l=12$ &  &  $AbAbCbAbdCbd$& &   $9^*2\,2$ \cr
  $l=10$  &  & $AbACdCbCdC$ &  & $.2 2:2$  &  &  $l=12$ &  & $AbAbCbCdCbCd$ & & $8^*2\,1\,1::2\,0$ \cr
    &  &  &  &   &  &  $l=12$ &  & $AbAbCbdCbCdC$ & &  $8^*2\,1\,1\,0:.2\,0$ \cr
  $l=11$  &  & $AbAbACbCdCd$ &  & $.2\,1\,1\,1\,1$  &  & $l=12$ &  &  $AbAbCbdCbdCd$& & $9^*2\,1\,1$  \cr
  $l=11$  &  & $AbAbCbCbdCd$ &  & $.2\,1\,1.2\,1\,0$  &  & $l=12$&  & $AbAbCdCbCdCd$ & &  $8^*2\,1\,1\,1\,0$  \cr
  $l=11$  &  & $AbAbCbdCbdC$ &  & $.2\,1\,1:2\,1$  &  &   $l=12$  &  &  $AbACbAdCbdCd$& &  $12$L \cr
  $l=11$  &  & $AbAbCdCbCdC$ &  & $.2\,1\,1\,1:2$  &  &  $l=12$ &  & $AbACbCbCbdCd$ & &  $8^*2\,1\,0.2\,1\,0$ \cr
  $l=11$  &  & $AbACbACbdCd$ &  & $9^*2\,1\,0$  &  &  $l=12$ &  &  $AbACbCbdCbCd$& & $9^*.2\,1:.2$  \cr
  $l=11$  &  & $AbACbCdCbCd$ &  & $8^*2\,1\,0::2\,0$  &  & $l=12$  &  & $AbACbCbdCbdC$ & & $8^*2\,1\,0:.2\,1\,0$  \cr
  $l=11$  &  & $AbACbCdCdCd$ &  & $.2\,2\,1\,1$  &  &  $l=12$ &  & $AbACbCdCbCdC$ & & $9^*2\,1:2$  \cr
  $l=11$  &  & $AbACbdCbCdC$ &  & $8^*2\,1:.2\,0$  &  &  $l=12$ &  & $AbACbCdCbdCd$ & & $10^{**}:2\,1\,0$  \cr
  $l=11$  &  & $AbACdCbCdCd$ &  & $8^*2\,2\,0$  &  &  $l=12$ &  & $AbACbdCbCdCd$ & & $10^{**}.2\,1$  \cr
    &  &  &  &   &  &  $l=12$ &  & $AbACbdCbdCdC$ & & $10^{**}:2\,1$  \cr
  $l=12$  &  & $AbAbAbACbdCd$ &  & $8^*2\,2\,1\,0$  &  &$l=12$ &  &  $AbCbAbCdCbCd$& &  $10^{**}:2 0::.2\,0$ \cr
  $l=12$  &  & $AbAbACbAbdCd$ &  & $9^*.2\,2$  &  & $l=12$  &  & $AbCbACbdCbCd$ & &  $10^{**}2\,0::.2\,0$ \cr}

\normalsize

\medskip

For $W=(Ab)^n$ ($n\ge 2$), $w_1=ACbdCdCd$ the family of basic
polyhedra beginning with $12$J ($AbAbACbdCdCd$) is obtained, and
for $W=(Ab)^n$ ($n\ge 1$), $w_1=ACbAdCbdCd$  the family of basic
polyhedra beginning with $12$L ($AbACbAdCbdCd$) is obtained.

\vspace*{1pt}\textlineskip  %) USE THIS MEASUREMENT WHEN THERE IS
\section{Applications of Minimum Braids and $BFR$s}    %) A SECTION HEADING
\vspace*{-0.5pt}

\vspace*{1pt}\textlineskip  %) USE THIS MEASUREMENT WHEN THERE IS
\subsection{Graph Trees}    %) A SECTION HEADING
\vspace*{-0.5pt}

A rational $KL$ in Conway notation is any sequence of natural
numbers not beginning or ending with 1, where each sequence is
identified with its inverse. From this definition is computed the
number of rational $KL$s with $n$ crossings. It is given by the
formula

$$2^{n-4}+2^{[n/2]-2}$$

\noindent that holds for every $n\ge 4$. This very simple formula
is derived first by C.~Ernst and D.W.~Sumners  in another form
[8], and later independently by S.~Jablan [9,10]. For $n\ge 4$ we
can compute the first 20 numbers of this sequence. After
prepending to it the first three numbers 1 for $n=1,2,3$, the
result is the sequence: 1, 1, 1, 2, 3, 6, 10, 20, 36, 72, 136,
272, 528, 1056, 1080, 4160, 8256, 16512, 32986, 65792, 131328,
262656, 524800, $\ldots$ This sequence is included in {\it On-Line
Encyclopedia of Integer Sequences}
(http://www.research.att.com/$^\sim $njas/sequences/) as the
sequence $A005418$. The number of rational knots with $n$
crossings ($n\ge 3$) is given by the formula

$${2^{n - 3} +
        2^{[{n\over 2}] -
              2^{{{(n - 1)}\pmod 2}}} + (-1)^{(n -
                  1){[{n\over 2}]\pmod 2}}}\over 3$$

\noindent so we can simply derive the formula for the number of
rational links with $n$ crossings as well.

A graph-theoretical approach to knot theory is proposed by
A.~Caudron [3]. T.~Gittings established a mapping between minimum
braids with $s$ strands and trees with $s+1$ vertices and
conjectured that the number of graph trees of $n$ vertices with
alternating minimum braids is equal to the number of rational
$KL$s with $n$ crossings [5, Conjecture 1].

\vspace*{1pt}\textlineskip  %) USE THIS MEASUREMENT WHEN THERE IS
\subsection{Amphicheiral $KL$s}    %) A SECTION HEADING
\vspace*{-0.5pt}

$KL$ is {\it achiral} (or {\it amphicheiral}) if its "left" and
"right" forms are equivalent, meaning that one can be transformed
to the other by an ambient isotopy.  If an oriented knot or link
$L$ can be represented by an antisymmetrical vertex-bicolored
graph on a sphere, whose vertices with the sign $+1$ are white,
and vertices with the sign $-1$ are black, it is achiral. In this
case, for an oriented knot or link $L$ there exists an
antisymmetry (sign-changing symmetry) switching orientations of
vertices, i.e., mutually exchanging vertices with the signs $+1$
and $-1$ [9,10]. In the language of braid words, this means that
its corresponding braid word is antisymmetric (or palindromic):
there exist a mirror antisymmetry transforming one letter to
another and {\it vice versa} and changing their case (i.e.,
transforming capital to lower case letters and {\it vice versa}).
For example, the reduced braid words $Ab\,|Ab$ or $ABac\,|BDcd$
are palindromic, where the anti-mirror is denoted by $|$. Hence,
we believe that the origin of all oriented achiral $KL$s are
palindromic reduced braids.

\medskip

{\bf Conjecture} An oriented $KL$ is achiral {\it iff} it can be
obtained from a palindromic reduced braid by a symmetric assigning
of degrees.

\medskip

For $s=2$ all alternating $BFR$s are of the form $(Ab)^n$ ($n\ge
2$), defining a series of the basic polyhedra $(2n)^*$, beginning
with $2\,2$, $.1=6^*$, $8^*$, $10^*$, $12^*$, {\it etc.} All of
them are achiral $KL$s, representing a source of other achiral
$KL$s. From 4:1-01 $AbAb$ ($2\,2$ or $4_1$) by a symmetric
assigning of degrees we can derive achiral alternating knots with
$n\le 10$ crossings: 6:1-02 $A^2bAb^2$ ($2\,1\,1\,2$ or $6_3$),
8:1-05 $A^3bAb^3$ ($3\,1\,1\,3$ or $8_9$), 10:1-017 $A^3b^2A^2b^3$
($(3,2)\,(3,2)$ or $10_{79}$), and one achiral alternating link
with $n\le 9$ crossings: 8:3-05a $A^2b^2A^2b^2$ (($2,2)\,(2,2)$ or
$8_4^3$), {\it etc.} In general, from $AbAb$ the following
families of achiral alternating $KL$s are derived:

\medskip

\halign{ # \hfill &  #  &  # \hfill  & # $\quad \quad $ & # \hfill
& # & # \hfill \cr
 $A^pbAb^p$ & & $p\,1\,1\,p$ & & $A^pb^qA^qb^p$ & &$(p,q)\,(p,q)$ \cr}

\medskip

Borromean rings 6:3-02 $AbAbAb$ ($.1=6^*$ or $6_2^3$) are the
origin of achiral alternating knots 8:1-07 $A^2bAbAb^2$ ($.2.2$ or
$8_{17}$), 10:1-020 $A^2bA^2b^2Ab^2$ ($.2.2.2\,0.2\,0$ or
$10_{99}$), 10:1-022 $A^2b^2AbA^2b^2$ ($2.2.2.2$ or $10_{109}$),
and of the link 8:3-04a $Ab^2AbA^2b$ ($.2:2\,0$ or $8_6^3$), {\it
etc.} In general, from $AbAbAb$ the following families of achiral
alternating $KL$s are derived:

\medskip

\halign{ # \hfill &  #  &  # \hfill  & # $\quad \quad $ & # \hfill
& # & # \hfill \cr
 $A^pbAbAb^p$ & & $.p.p$ & & $A^pbA^qb^qAb^p$ & &$.p.p.q\,0.q\,0$ \cr
 $Ab^pAbA^pb$ & & $.p:p\,0$ & & $A^pb^qA^rb^rA^qb^p$ & &$p.q.r.r.q.p$ \cr
 $A^pb^qAbA^qb^p$ & & $p.q.q.p$ & & $$ & &$$ \cr}

\medskip

Achiral basic polyhedron $AbAbAbAb$ ($8^*$) is the origin of the
following families of alternating achiral $KL$s:

\medskip

\halign{ # \hfill &  #  & # \hfill  & # $\quad \quad $ & # \hfill
& # & # \hfill \cr
 $A^pbAbAbAb^p$ & & $8^*p.p$ & & $A^pbAb^qA^qbAb^p$ & &$8^*p.q:.q.p$ \cr
 $AbA^pbAb^pAb$ & & $8^*p:.p$ & & $A^pb^qA^rbAb^rA^qb^p$ & &$8^*p.q.r.r.q.p$ \cr
 $A^pb^qAbAbA^qb^p$ & & $8^*p.q.q.p$ & & $A^pbA^qb^rA^rb^qAb^p$ & &$8^*p.q.q.p:r.r$ \cr
 $A^pbA^qbAb^qAb^p$ & & $8^*.p:q.q:p$ & & $$ & &$$ \cr}

\medskip

In the same way it is possible to derive achiral alternating $KL$s
from all achiral basic polyhedra $(Ab)^n$ for $n\ge 5$.

From the antisymmetry condition it follows that every palindromic
braid has an even number of strands. For $s=4$ and $l\le 12$
palindromic algebraic generating braids are:

$AbACbdCd$, $l=8$ with the corresponding achiral link
$2\,2\,2\,2$,

$AbAbCbCdCd$, $l=10$, with the corresponding achiral knot
$2\,1\,1\,1\,1\,1 \,1\,2$.

\noindent The palindromic polyhedral generating braids are:

$AbACbCbdCd$, $l=10$, with the corresponding achiral knot
$.2\,1.2\,1$,

$AbAbACbdCdCd$, $l=12$, with the corresponding achiral link $12$J,

$AbACbAdCbdCd$, $l=12$, with the corresponding achiral knot $12$L,

$AbACbCbCbdCd$, $l=12$, with the corresponding achiral link
$8^*2\,1\,0.2\,1\,0$,

$AbCbAbCdCbCd$, $l=12$, with the corresponding achiral knot
$10^{**}:2 0::.2\,0$,

$AbCbACbdCbCd$, $l=12$, with the corresponding achiral knot
$10^{**}2\,0::.2\,0$.

From the generating braid $AbACbdCd$ following families of
alternating achiral $KL$s are derived:

\medskip

\footnotesize

\halign{ # \hfill &  #  & # \hfill  & #  & # \hfill & # & # \hfill
\cr
 $A^pbACbdCd^p$ & & $p\,1\,2\,2\,1\,p$& & $A^pbA^qC^rb^rd^qCd^p$ & &$(((p,q),r)+)\,(((p,q),r)+)$ \cr
 $AbAC^pb^pdCd$ & & $(p,2+)\,(p,2+)$ & & $A^pb^qAC^rb^rdC^qd^p$ & &$(q,p\,1,r)\,(q,p\,1,r)$ \cr
 $Ab^pAC^qb^qdC^pd$ & & $(p,q,2)\,(p,q,2)$ & & $A^pb^qA^rC^sb^sd^rC^qd^p$ & &$(q,(p,r),s)\,(q,(p,r),s)$ \cr}

\normalsize

\medskip

From the same palindromic non-alternating generating braid the
following families of achiral $KL$s are obtained:

\medskip

\halign{ # \hfill &  #  & # \hfill  & # $\quad \quad $ & # \hfill
& # & # \hfill \cr
 $A^pBacBDcd^p$ & & $p\,p$& & $A^pBacBDcd^p$ & &$2\,p\,p\,2$ \cr
 $AbAc^pB^pdCd$ & & $(p,2)\,(q,2)$ & & $A^pbAc^qB^qdCd^p$ & &$(p\,1,q)\,(p\,1,q)$ \cr}

\medskip

In the same way is possible to continue the derivation of achiral
$KL$s from other palindromic reduced braids.

The family of achiral odd crossing number knots discovered by
J.~Hoste, M.~Thistlethwaite and J.~Weeks in 1998 [11] can be
extended to the two-parameter $BFR$ defined by the palindromic
braid $ABaB^qC^pBAdcb^pc^qDcd$ corresponding to the family of
non-alternating achiral odd-crossing knots with $n=7+4p+4q$
crossings
$$10^{**}(-2p)\,0.-1.-2\,0.(2q):(-2p)\,0.- 1.-2\,0.(2q)$$.

\vspace*{1pt}\textlineskip  %) USE THIS MEASUREMENT WHEN THERE IS
\subsection{Unlinking Numbers and Unlinking Gap}    %) A SECTION HEADING
\vspace*{-0.5pt}

T.~Gittings [5] noticed that it might be possible to calculate
unlinking numbers from minimum braids. Unfortunately, this is true
only for $KL$s with $n\le 10$ crossings, including the link
$4\,1\,4$ ($9_4^2$) and the Nakanishi-Bleiler example $5\,1\,4$
($10_8$) with an unlinking gap [12].

\medskip

\begin{definition}
The {\it minimum braid unlinking gap} is the positive difference
between the unlinking number obtained from a minimum braid
$u_{B}(L)$ and unlinking number $u(L)$ of a link $L$, i.e.,

$$\delta _B=u_{B}(L)-u(L)>0.$$
\end{definition}

\medskip

The unlinking gap [12] for minimum braids appears for $n=11$. The
following alternating links given in Conway notation, followed by
their minimum braids have the minimum braid unlinking gap:

\medskip

\halign{ # \hfill &  # & # &  # \hfill  & $\quad \quad \quad $# &
# \hfill &  # & # \hfill  \cr
  $.5.2$ & & $A^5bAbAb^2$& & &$8^*3.2$ & &$A^3bAbAbAb^2$ \cr
  $.3.4$& & $A^4bAbAb^3$& & & $8^*3:2$& & $A^3bA^2bAbAb$\cr
  $8^*4$ & & $A^4bAbAbAb$& & &$8^*2.2:.2$ & &$A^2bAbA^2bAb^2$ \cr
  $.2.3.3\,0$ & &$A^3bA^3bAb^2$ & & &$10^*2$ & &$A^2bAbAbAbAb$ \cr}

\medskip

\noindent For the links $.5.2$, $.3.4$ the value of minimum braid
unlinking gap is $\delta _B=2$, and for other links from this list
$\delta _B=1$. Hence, minimum braid unlinking number is different
from the unlinking number and represents a new $KL$ invariant.

\vspace*{1pt}\textlineskip  %) USE THIS MEASUREMENT WHEN THERE IS
\subsection{Periodic Tables of $KL$s}    %) A SECTION HEADING
\vspace*{-0.5pt}

Periodic tables of $KL$s can be established in three ways:
starting with families of $KL$s given in Conway notation
[9,10,13], with minimum braids [5], or with $BFR$s. Since we have
established correspondence between $BFR$s and $KL$s in Conway
notation, it follows that the same patterns (with regard to all
$KL$ polynomial invariants and $KL$ properties) will appear in all
cases. For example, for every family of $KL$s is possible to
obtain a general formula for Alexander polynomials, with
coefficients expressed by numbers denoting tangles in Conway
symbols, or from their corresponding parameters from minimum
braids or from $BFR$s. The same holds not only for $KL$
polynomials, but for all other properties of $KL$s: writhe,
amphicheirality, number of projections, unlinking number,
signature, periods, {\it etc.} [9,10,13].

\nonumsection{Acknowledgements}

We would like to express our gratitude to Thomas Gittings for his
critical reading of the manuscript, corrections, advice and
suggestions.

\nonumsection{References}

\end{document}